\documentclass[11pt]{article}

\usepackage{amsmath,amsthm,amsfonts,amscd,bezier}
\usepackage{amssymb}
\usepackage[english]{babel}
\usepackage[T1]{fontenc}
\usepackage{comment}
\usepackage[usenames]{color}
\usepackage{enumerate}

\newtheorem{teo}{Theorem}[section]

\newtheorem{lema}[teo]{Lemma}
\newtheorem{cor}[teo]{Corollary}
\newtheorem{prop}[teo]{Proposition}
\newtheorem{ex}[teo]{Example}
\newtheorem{rem}[teo]{Remark}
\newtheorem{conj}[teo]{Conjecture}
\newcommand{\pf}{\noindent{\bf Proof:\ \ }}

\newcommand{\fim}{\hfill $\rule{2.5mm}{2.5mm}$}
\DeclareMathOperator{\Oh}{{\mathcal O}}
\DeclareMathOperator{\pr}{pr}
\DeclareMathOperator{\I}{{\mathcal I}}
\DeclareMathOperator{\T}{{\mathcal T}}
\DeclareMathOperator{\J}{{\mathcal J}}
\DeclareMathOperator{\K}{{\mathcal K}}
\DeclareMathOperator{\N}{{\mathcal N}}
\DeclareMathOperator{\Ss}{{\mathcal S}}
\DeclareMathOperator{\F}{{\mathcal F}}
\begin{document}

\title{Colengths of fractional ideals and Tjurina number of a reducible plane curve}
\author{Abramo Hefez and Marcelo Escudeiro Hernandes}
\maketitle

\begin{abstract}
In this work, we refine a formula for the Tjurina number of a reducible algebroid plane curve defined over $\mathbb C$ obtained in the more general case of complete intersection curves in \cite{Bayer}. As a byproduct, we answer the affirmative to a conjecture proposed by A. Dimca in \cite{Dimca}. Our results are obtained by establishing more manageable formulas to compute the colengths of fractional ideals of the local ring associated with the algebroid (not necessarily a complete intersection) curve with several branches.  We then apply these results to the Jacobian ideal of a plane curve over $\mathbb C$ to get a new formula for its Tjurina number and a proof of Dimca's conjecture. We end the paper by establishing a connection between the module of K\"ahler differentials on the curve modulo its torsion, seen as a fractional ideal, and its Jacobian ideal, explaining the relation between the present approach and that of \cite{Bayer}.  
\end{abstract}

\section{Introduction}

Many invariants of an algebroid curve $C$, defined over an algebraically closed field $\mathbb K$, may be computed using colengths of fractional ideals of its local ring $\mathcal{O}$. This is so for the Tjurina number of a complete intersection curve $C$ when $\mathbb K=\mathbb C$. For a reduced curve $C$ with irreducible components $C_i,\ 1\leq i\leq r$, it is natural to seek for relations among the invariants of $C$ and $C_i$ as in the emblematic example  for the Milnor number $\mu(C)$ of a plane curve $C$:  $$\mu(C)=\sum_{i=1}^{r}\mu(C_{i})+2\sum_{1\leq i<j\leq r}\textup{I}(C_i,C_j)-r+1,$$ where $\textup{I}(C_i,C_j)$ is the intersection multiplicity of the branches $C_i$ and $C_j$ (see \cite[Theorem 6.5.1]{wall}, for instance).

The value set $E$ of a fractional ideal $\mathcal{I}$ is an important object that can be used to compute the colength of $\mathcal{I}$ in $\mathcal{O}$, as shown in \cite{danna}, \cite{D87} and \cite{GH}, for instance. The description of $E$ can be done via algorithms (see \cite{emilio2}), or using a recursive method that depends on particular elements in $E$ (the relative maximal points) as described  in \cite{D87} and \cite{GH}. 

In Section 2, we present characteristic free general results concerning value sets of fractional ideals of the local ring of any curve that will allow us to give an alternative method to compute their colength, without being necessary to know the relative maximal points of their value sets in advance (see Corollary \ref{theorem-colenth}).

In Section 3, by considering a plane algebroid curve $C$, defined over $\mathbb C$, with irreducible components $C_i,\ 1\leq i\leq r$,  we specialize the above mentioned results to the value set of the Jacobian ideal to rephrase a formula presented in \cite{Bayer} that relates the Tjurina number $\tau(C)$ of $C$ and the Tjurina number $\tau(C_i)$ of $C_i$ (see Proposition \ref{tjurinaref} and Theorem \ref{tjurinaIJK}). In particular, this allows us to prove in Section 4 the  following theorem (cf. Theorem \ref{proof-dimca}) that was conjectured by A. Dimca in \cite[Conjecture 1.7]{Dimca}:

\smallskip

\noindent{\sc\bf  Theorem} {\it Let $J\cup K$ be a partition of $\{1,\ldots,r\}$. If $C_J=\bigcup_{j\in J}C_j$ and $C_K=\bigcup_{k\in K}C_k$,  then} $$\tau(C_J\cup C_K)\leq\tau(C_J)+\tau(C_K)+2\cdot \textup{I}(C_J,C_K)-1,$$
{\it where $\textup{I}(C_J,C_K)$ denotes the intersection multiplicity of $C_J$ and $C_K$. In particular, }
$$\tau(C)-\sum_{1\leq i\leq r}\tau(C_{i})\leq 2\sum_{1\leq i<j\leq r}\textup{I}(C_i,C_j)-r+1=\mu(C)-\sum_{1\leq i\leq r}\mu(C_{i}).$$
\smallskip

In Section 5, we build a bridge between the approach of \cite{Bayer} that used K\"ahler differentials, and the present, where we use the Jacobian ideal, to get the formula that relates Tjurina's number of a plane curve and the Tjurina numbers of its components.

\section{Fractional Ideals of local rings of curves}

In this section, we keep the same framework and notations as in \cite{Bayer}. Let $C$ be an algebroid reduced curve defined over an algebraically closed field $\mathbb{K}$ in an $n$-dimensional Euclidean space defined by a radical ideal ${\mathfrak A}\subset\mathbb{K}[[X_1,\ldots ,X_n]]$. Here, all arguments and results are characteristic-free.

Let $\wp_1,\ldots,\wp_r$ be the minimal primes of $\mathcal{O}=\mathbb{K}[[X_1,\ldots ,X_n]]/\mathfrak A$ and set $I=\{1,\ldots,r\}$. For any non-empty subset $J$ of $I$, put $\mathcal{O}_J=\mathcal{O}/\cap_{j\in J} \wp_{j}$. We denote the local ring $\mathcal{O}_{\{i\}}$ associated by the branch $C_i$ defined by $\wp_i$ simply by $\mathcal{O}_i$. Notice that $\mathcal{O}_I=\mathcal{O}$. 

If $\widetilde{\mathcal{O}_J}$ is the integral closure of $\mathcal{O}_J$ in its total ring of fractions $\mathcal{K}_J$, one has $\widetilde{\mathcal{O}_J}\simeq \prod_{j\in J} \widetilde{\mathcal{O}_{j}} \simeq \prod_{j\in J} \mathbb{K}[[t_{j}]]$ and $\mathcal{K}_J \simeq \prod_{j\in J} \mathbb{K}((t_{j}))$. 

We put $\mathcal{K}=\mathcal{K}_I=\prod_{i\in I} \mathbb{K}((t_{i}))$ and we denote by $\pi_J:\mathcal{K}\rightarrow\mathcal{K}_J$ the natural projection, where $\pi_i$ stands for $\pi_{\{i\}}$. Notice that the regular elements (i.e., the non-zero divisors) in $\mathcal K$ are the $r$-tuples with none of their components being zero.

Considering the natural surjections $\Oh \to \Oh_i$  and viewing $\mathcal{O}_i\subseteq\widetilde{\mathcal{O}_i}\subset\mathcal{K}_i$, we get   a ring monomorphism 
\[ \mathcal{O} \hookrightarrow \prod_{i\in I} \Oh_i \hookrightarrow \prod_{i\in I} \widetilde{\Oh_i} \hookrightarrow \K \] that allows us to identify $\mathcal{O}$ with a subring of $\mathcal{K}\simeq \mathbb K((t_1))\times \cdots \times \mathbb K((t_r))$, which makes things much easier, because of the simple structure of $\K$. With this identification, for $J\subset I$, the natural surjection $\Oh \to \Oh_J$ corresponds to the restriction $\pi_J\colon \Oh \to \Oh_J$ of the projection $\pi_J\colon \K \to \K_J$. 

A fractional ideal $\mathcal{I}$  of $\mathcal{O}$  is an $\mathcal{O}$-submodule of $\K$ that contains a regular element $g\in\mathcal{O}$ satisfying $g\mathcal{I}\subseteq\mathcal{O}$. For such $\I$ we define $\I_J\colon= \pi_J(\I)\subset \K_J$. Observe that $\I_J$ is a fractional ideal of $\mathcal{O}_J$.

Let $\overline{\mathbb{Z}}:=\mathbb{Z}\cup\{\infty\}$ and let $\nu_i:\mathcal{K}_i\rightarrow\overline{\mathbb{Z}}$ be the discrete valuation given by $\nu_i\left ( \frac{g(t_i)}{h(t_i)}\right )=ord_{t_i}(g(t_i))-ord_{t_i}(h(t_i))$, where $ord_{t_i}(g(t_i))=\infty$ if and only if $g(t_i)=0$. 
We endow the Cartesian power $\overline{\mathbb{Z}}^s$ with the product order. 

Given $J=\{j_1<\cdots <j_s\}\subseteq I$  we define
$$\begin{array}{lllc}
	\nu_J: & \mathcal{K}_J & \rightarrow & \overline{\mathbb{Z}}^s \\
	& z & \mapsto & (\nu_{j_1}(\pi_{j_1}(z)),\ldots ,\nu_{j_s}(\pi_{j_s}(z)))
\end{array}$$
and $E_J:=\nu_J(\mathcal{I}_J)$ is called the value set of the factional ideal $\mathcal{I}_J$ of ${\mathcal O}_J$. We denote $\nu_I$ simply by $\nu$,  $\nu(\mathcal{I})$ by $E$ and for $J=\{i\}$ we write $\nu_{\{i\}}=\nu_i$ and $E_{\{i\}}=E_i$. In particular, $\Gamma:=\nu(\mathcal{O})$ is the value semigroup of $C$ and $\Gamma_i=\nu_i(\mathcal{O}_i)$ is the corresponding value semigroup associated to the branch $C_i$. 

If $r>1$, the semigroup $\Gamma$ is not finitely generated, but considering
$$\inf(\alpha,\beta)=(\min\{\alpha_1,\beta_1\},\ldots ,\min\{\alpha_r,\beta_r\})$$
for $\alpha=(\alpha_1,\ldots ,\alpha_r), \beta=(\beta_1,\ldots ,\beta_r)\in\Gamma$ then $(\Gamma,\inf ,+)$ is a finitely generated semiring (see \cite{emilio1}) and, since $\mathcal{I}$ is an $\mathcal{O}$-module, it follows that $(E,\inf,+)$ is a (finitely generated) $\Gamma$-semimodule, that is, $\gamma+\delta\in E$ for every $\gamma\in \Gamma$ and $\delta\in E$ (see \cite{emilio2}).

In order to simplify the notation we put $\nu_i(z):=\nu_i(\pi_i(z))$, where $z\in\mathcal{I}$.

In what follows, we present some properties concerning the value set of a fractional ideal that will be used in the sequel, referring the reader, for more details, to \cite{danna}, \cite{D87} and \cite{GH}.

There exists a unique element $\min(E)=\delta^0\in E$ such that $\delta-\delta^0\in\mathbb{N}^r$ for every $\delta\in E$ called {\it minimum of} $E$. There exists an element $\beta\in \mathbb Z^r$ such that $\beta+\mathbb{N}^r\subseteq E$, hence $\beta \in E$ and the minimum of such $\beta$ is called the {\it conductor} of $E$ and denoted by $c(E)$. In particular, if $\kappa\in\mathcal{K}$ with $\nu(\kappa)\geq c(E)$ then $\kappa\in \mathcal I$.

Given $\emptyset\neq J=\{j_1<\cdots <j_s\}\subseteq I$ if $\pr_J:\overline{\mathbb{Z}}^r\rightarrow\overline{\mathbb{Z}}^s$ denotes the projection $\pr_J(\alpha_1,\ldots ,\alpha_r)=(\alpha_{j_1},\ldots ,\alpha_{j_s})$ then one has
\begin{equation}\label{compartative}
	\begin{array}{ll}
E_J=\pr_J(E),\ \ & \min(E_J)=\min(\pr_J(E)),\vspace{0.2cm}\\
E\subset\prod_{i\in I} E_i,\ \ & c(E_J)\leq \pr_J(c(E)).
\end{array}\end{equation}
In addition, for $\alpha\in\mathbb{Z}^r$,  we define the $J$-fiber of $\alpha$ in $E$ as
$$F_J(E,\alpha)=\{\beta\in E;\ \pr_J(\beta)=\pr_J(\alpha)\ \mbox{and}\ \pr_{I\setminus J}(\beta)>\pr_{I\setminus J}(\alpha)\}.$$

We say that $\alpha\in E$ is a {\it maximal point} of $E$ if $\bigcup_{i\in I}F_{\{i\}}(E,\alpha)=\emptyset$. A maximal point $\alpha\in E$ is a {\it relative maximal} if $F_{J}(E,\alpha)\neq\emptyset$ for every $J\subseteq I$ with $\sharp J\geq 2$. We denote by $M(E)$ and $RM(E)$ the sets of maximal and relative maximal of $E$, respectively. In particular, if $r=2$ then $M(E)=RM(E)$. 

If $\min(E)=(\delta^0_{1},\ldots ,\delta^0_{r})$ and $c(E)=(\delta^c_{1},\ldots ,\delta^c_{r})$ denote the minimum and the conductor of $E$, respectively, then 
$$RM(E)\subset M(E)\subset [\delta^0_{1},\delta^c_{1}]\times\cdots\times [\delta^0_{r},\delta^c_{r}],$$ 
which implies that $E$ has  finitely many maximal points.

If $\J,\I\subseteq \K$ are fractional ideals of $\mathcal{O}$ with $\mathcal{J}\subseteq \mathcal{I}$,  $D=\nu(\mathcal{J})$ and $E=\nu(\mathcal{I})$, then the length (as $\mathcal{O}$-module) $\ell_{\mathcal{O}}\left (\frac{\mathcal{I}}{\mathcal{J}}\right )$ can be computed by means $RM(D_J), RM(E_J)$, $D_i$ and $E_i$ where $\emptyset\neq J=\{j_1,\ldots ,j_s\}\subseteq I,\ 1<s\leq r$ and $i\in I$. 

Indeed, defining $\pr_*(\alpha_1,\ldots ,\alpha_s)=\alpha_s$, and 
\begin{equation} \label{thetai}
\Theta_1(E)=0\ \ \ \mbox{and}\ \ \ \Theta_i(E)=\sharp\left ( \bigcup_{i\in J\subseteq \{1, \ldots , i\}}\pr_*(RM(E_J))\right ),\ \ 2\leq i\leq r,
\end{equation}
the following result was proved in \cite{GH}.\smallskip

\noindent \textbf{Theorem} (\cite[Corollary 11]{GH}) \emph{Let $\I$ be a fractional ideal of $\Oh$ with value set $E$, then for every $\gamma \geq c(E)$ one has
\[
\ell_\mathcal{O}\left(\dfrac{\mathcal I}{\mathcal I(\gamma)} \right)=\sum_{i=1}^r (\gamma_i-\delta_i^0 -\sharp {\mathcal G}(E_i) -\Theta_i(E)),
\]
where $(\delta_1^0, \ldots,\delta_r^0):=\min(E)$, ${\mathcal G}(E_i):=\{z\in\mathbb{Z};\ \delta^0_{i}\leq z\notin E_i\}$, the set of {\it gaps} of $E_i$, and $\I(\gamma)=\{z\in\mathcal{I};\ \nu(z)\geq \gamma\}$. }\medskip

\begin{rem} \label{r=1} Let $\J \subset \I$ be $\mathbb K$-vector subspaces of $\mathbb K((t))$, then it is well known and easy to prove that
\[
\dim_{\mathbb K} \dfrac{\I}{\J}=\sharp (\nu(\I) \setminus \nu(\J)).
\]
\end{rem}

Now, as an application of the above theorem, we get the following consequence:

\begin{prop}\label{codim}
Let  $\mathcal{J},\mathcal{I}\subseteq\mathcal{K}$ be fractional ideals of $\mathcal{O}$ with $\mathcal{J}\subseteq\mathcal{I}$,  $D=\nu(\mathcal{J})$ and $E=\nu(\mathcal{I})$, then 
$$
\ell_{\mathcal{O}}\left (\frac{\mathcal{I}}{\mathcal{J}}\right )  =\sum_{i\in I}\left (\ell_i\left ( \frac{\mathcal{I}_i}{\mathcal{J}_i}\right )+\Theta_i(D)-\Theta_i(E)\right ),
$$
where $\ell_i\left ( \frac{\mathcal{I}_i}{\mathcal{J}_i}\right )$ denotes the length of $\frac{\mathcal{I}_i}{\mathcal{J}_i}$ as $\mathcal{O}_i$-module.
\end{prop}

\pf
According to \cite[Proposition 2.2]{danna}, if $\mathcal{J}\subseteq\mathcal{I}$ are fractional ideals with value sets $D$ and $E$, respectively, then \begin{equation}\label{quotient}
\ell_{\mathcal{O}}\left ( \frac{\mathcal{I}}{\mathcal{J}}\right )=\ell_{\mathcal{O}}\left ( \frac{\mathcal{I}}{\mathcal{I}(\gamma)}\right )-\ell_{\mathcal{O}}\left ( \frac{\mathcal{J}}{\mathcal{J}(\gamma)}\right )
\end{equation}
for any $\gamma=(\gamma_1,\ldots ,\gamma_r)\geq c(D)\geq c(E)$.

By the above stated theorem, if $\min(E)=(\delta^0_{E,1},\ldots ,\delta^0_{E,r})$, then
\begin{equation}\label{quotientI}
\ell_{\mathcal{O}}\left ( \frac{\mathcal{I}}{\mathcal{I}(\gamma)}\right )=\sum_{i\in I}\left (\gamma_i-\delta^0_{E,i}-\sharp {\mathcal G}(E_i)-\Theta_i(E)\right ).	
\end{equation}

Replacing (\ref{quotientI}) (and the corresponding expression for $\mathcal{J}$) in (\ref{quotient}), and observing that $\sharp(E_i\setminus D_i)=\delta^0_{D,i}-\delta^0_{E,i}+\sharp {\mathcal G}(D_i)-\sharp {\mathcal G}(E_i)$, we get
$$\ell_{\mathcal{O}}\left (\frac{\mathcal{I}}{\mathcal{J}}\right )  =\sum_{i\in I}\left (\sharp (E_i\setminus D_i)+\Theta_i(D)-\Theta_i(E)\right ).$$

Now, since by Remark \ref{r=1}, one has that $\ell_i\left ( \frac{\mathcal{I}_i}{\mathcal{J}_i}\right )=\dim_{\mathbb K}\frac{\I_i}{\J_i}=\sharp (E_i\setminus D_i)$, we are done.

\fim

\begin{ex}\label{prod} 	Let $\mathcal{I}$ be a fractional ideal of $\mathcal{O}$ with value set $E\subseteq \overline{\mathbb{Z}^r}$. Since $\nu\left (\prod_{i\in I}\mathcal{I}_i\right )=\prod_{i\in I}E_i$ has no maximal points, and since $\left (\prod_{i\in I}\mathcal{I}_i\right )_j=\mathcal{I}_j$ and $\Theta_1(E)=0$, then, by Proposition \ref{codim}, we get
$$\ell_{\mathcal{O}}\left (\frac{\prod_{i\in I}\mathcal{I}_i}{\mathcal{I}}\right ) =\sum_{i=1}^{r}\Theta_i(E)=\sum_{i=2}^{r}\Theta_i(E).	
$$	
\end{ex}

Notice that, at this point, to compute $\Theta_i(D)$ and $\Theta_i(E)$, we need to determine the sets of relative maximals of all the projections of $D$ and $E$, which is not suitable for our purpose. So, in the sequel, we will give an alternative way to compute the numbers $\Theta_i(E)$ for the  value set $E$ of a fractional ideal.

Given $\alpha\in\mathbb{Z}^r$ and $i\in I$ we define the \emph{closed $i$-fiber} of $\alpha$ in $A\subseteq\overline{\mathbb{Z}}^r$ as  
$$\overline{F}_i(A,\alpha)=\{\beta\in A;\ \beta_i=\alpha_i\ \mbox{and}\ \beta_j\geq\alpha_j\ \mbox{for}\ j\neq i\}.$$

According to \cite[Lemma 9]{GH}, for $1<i\leq r$, one has
\begin{equation}\label{theta}
\begin{array}{ll}
	\Theta_i(E) & =\sharp\{ z\in E_i\ \mbox{and}\ \overline{F}_i(E_{\{1,\ldots ,i\}},(\delta^c_1,\ldots ,\delta^c_{i-1},z))=\emptyset\}\\ \\
	& =\sharp E_i\setminus \{ z<\delta_i^c;\ \overline{F}_i(E_{\{1,\ldots ,i\}},(\delta^c_1,\ldots ,\delta^c_{i-1},z))\neq\emptyset\},
\end{array}
\end{equation} 
where $(\delta^c_1,\ldots ,\delta^c_r)=c(E)$,

For the new description of the $\Theta_i(E)$, where $E$ is the value set of a fractional ideal $\mathcal I$,  that we have in mind, the 
zero-divisors 
in the fractional ideal will play an important role.

For a fractional ideal $\I$ of $\Oh$ (both viewed in $\mathcal K$) and for $i\in I$, we define
\[
\N_i(\I)=\I \cap \ker{\pi_i},
\]
where $\pi_i\colon \mathcal K \to \mathcal K_i$ is, as before, the projection on the $i$-th factor of $\K$. Now, for $J \subset I$ we define $\N_J(\I)=\bigcap_{i\in J}\N_i(\I)$.

Notice that $\N_J(\I)$ is a fractional ideal of $\Oh$. If in a discussion $\I$ is fixed, we denote $\N_J(\I)$ simply by $\N_J$. We also have that  $\nu_i(\N_J)\subset E_i$. In particular,
the value set $N_i:=\nu_{i}(\N_{I\setminus \{i\}}) \subset E_i$ has a conductor $c(N_i)=\rho_i$ which is such that $\rho_i\geq c(E_i)$.

\begin{lema}\label{condutor} With the above notation we have 
\begin{equation}\label{inclusions}
\prod_{i\in I} N_i \ \subset\ E\ \subset\ \prod_{i\in I}E_i.
\end{equation}
In particular, we get $(c(E_1),\ldots ,c(E_r))\leq c(E)\leq (\rho_1,\ldots ,\rho_r)$.	
\end{lema}
\pf  For each $i\in I$, take $\eta_i\in \N_{I\setminus \{i\}} \subset\mathcal{I}$, hence 
 $\sum_{i\in I}\eta_i\in\mathcal{I}$ and $(\nu_1(\eta_1),\ldots ,\nu_r(\eta_r))=\nu(\sum_{i\in I}\eta_i)\in E$, that is,
$$\prod_{i\in I} N_i=\prod_{i\in I}\nu_i(\N_{I\setminus \{i\}})\subset E.$$ 
The other inclusion follows from (\ref{compartative}).

The inequalities $(c(E_1),\ldots ,c(E_r))\leq c(E)\leq (\rho_1,\ldots ,\rho_r)$ are consequences of (\ref{inclusions}).
\fim \smallskip

The following result gives us more information about the value set $E$ by means of the value sets $N_i=\nu_i(\N_{{I\setminus \{i\}}})$, for $i\in I$.

\begin{prop}\label{box-out}
Let $\rho_j=c(N_j)$, for all $j\in I$. Given $\delta\in\mathbb{Z}$ and defining $\alpha=(\rho_1, \ldots, \rho_{i-1},\delta,\rho_{i+1}, \ldots,\rho_r)$, the following statements are equivalent:
\begin{enumerate}
	\item $\delta\in N_i$;
	\item $\overline{F}_i(\mathbb{Z}^r,\alpha) \subset E$;
	\item  $\overline{F}_i(\mathbb{Z}^r,\alpha) \cap E \neq \emptyset$;
	\item $\alpha\in E$.
\end{enumerate}
\end{prop}
\pf (1 $\Rightarrow$ 2) \, Since $\delta\in N_i$, and $\rho_j =c(N_j)$, then for all $n_j\in \mathbb N$, $j\in I\setminus \{i\}$, from \eqref{inclusions} one has
\[
(\rho_1+n_1,\ldots ,\rho_{i-1}+n_{i-1},\delta,\rho_{i+1}+n_{i+1},\ldots ,\rho_r+n_r)\in \prod_{i\in I} N_i \subset E.
\]
(2 $\Rightarrow$ 3) Immediate

\noindent (3 $\Rightarrow$ 4) Suppose that for some $n_1,\ldots,n_{i-1},n_{i+1},\ldots,n_r\in \mathbb N$ one has
\[ 
\beta=(\rho_1+n_1,\ldots ,\rho_{i-1}+n_{i-1},\delta,\rho_{i+1}+n_{i+1},\ldots ,\rho_r+n_r)\in E.
\]
Taking $\epsilon \in N_i$ such that $\epsilon \geq \delta$, we have, again from \eqref{inclusions}, that
\[
\gamma =(\rho_1,\ldots ,\rho_{i-1},\epsilon,\rho_{i+1},\ldots ,\rho_r) \in \prod_{i\in I}^r N_i \subset E.
\]
It follows that
\[
\alpha= (\rho_1,\ldots ,\rho_{i-1},\delta,\rho_{i+1},\ldots ,\rho_r)=\min\{\beta,\gamma\} \in E.
\]
\noindent (4 $\Rightarrow$ 1) Suppose that $\alpha=(\rho_1,\ldots ,\rho_{i-1},\delta,\rho_{i+1},\ldots ,\rho_r)\in E$ and let $z\in \mathcal I$ be such that $\nu(z)=\alpha$. Since, for $j\neq i$, we have $\nu_j(z)=\rho_j=c(N_j)$, there exists $\eta_j\in \N_j$ such that $\pi_j(z)=\pi_j(\eta_j)$. It follows that $z-\eta_j\in \N_j$, for every $j\in I\setminus \{i\}$, hence $y:=z-\sum_{j\in I\setminus \{i\}} \eta_j \in \N_{I\setminus \{i\}}$, and 
\[
\delta=\nu_i(z)=\nu_i(y)\in N_i.
\]
\fim

The above proposition leads to the following:

\begin{cor} \label{gendanna}
We have that $c(E)=(\rho_1,\ldots ,\rho_r)$, where $\rho_i=c(N_i)$. 
\end{cor}
\pf By Lemma \ref{condutor} we get $c(E)\leq (\rho_1,\ldots ,\rho_r)$.

Now, since $\rho_i-1\not\in N_i$, by Proposition \ref{box-out}, we have $(\rho_1,\ldots ,\rho_{i-1},\rho_i-1,\rho_{i+1},\ldots ,\rho_r)\not\in E$ for every $i\in I$. This implies that $(\rho_1,\ldots ,\rho_r) \leq c(E)$. 

Hence, $c(E)=(\rho_1,\ldots ,\rho_r)$.
\fim

\begin{rem} \label{Danna}If we apply Corollary \ref{gendanna} to the case $\mathcal{I}=\mathcal{O}$, we have that $E=\Gamma$ and $\N_j=\wp_j$, where $\wp_1,\ldots ,\wp_r$ denote the minimal primes of $\mathcal{O}$. In this case, $\N_{I\setminus \{i\}}=\bigcap_{j\in I\setminus \{i\}}\wp_j$, and, therefore,
$$
c(\Gamma)=\left ( c\left (\nu_1\left ( \cap_{j\in I\setminus\{1\}}\wp_j\right )\right ),\ldots ,c\left (\nu_r\left ( \cap_{j\in I\setminus\{r\}}\wp_j\right )\right )\right ),$$
as showed by D'Anna in \cite[Proposition 1.3]{danna}. In particular, one gets
\begin{equation} \label{conds}
c(\Gamma)_i= \nu_i\left ( \cap_{j\in I\setminus\{i\}}\wp_j\right )+\mu_i,
\end{equation}
where $\mu_i$ is the conductor of $\Gamma_i=\nu_i(\Oh_i)$.

\end{rem}

Notice that our Corollary \ref{gendanna} generalizes the content of the above remark to any fractional ideal $\mathcal I$ of $\mathcal O$.

For $a,b\in I=\{1,\ldots,r\}$ with $a<b$, we denote by $[a,b)$ the semi-open interval in $I$ with extremes $a$ and $b$. With the notations introduced in this section, we get the following rephrasing of the $\Theta_i(E)$, defined in \eqref{thetai} in a combinatorial way, into algebraic expressions which are more convenient to our purpose.

\begin{cor}\label{thetaf} For $1< i\leq r$ we have
$$\Theta_i(E) =\sharp \left(E_i\setminus  \nu_i\left ( \N_{[1,i)} \right ) \right) \ \ \left(=\ell_i \left(\frac{\I_i}{\pi_i(\N_{[1,i)})}\right)\right).$$
\end{cor}
\pf From \eqref{theta} and the equivalence in Proposition \ref{box-out}, we have that
$$
	\begin{array}{ll}
		\Theta_i(E) & =\sharp \left(E_i\setminus \{ z<\delta_i^c;\ \overline{F}_i(E_{\{1,\ldots ,i\}},(\delta^c_1,\ldots ,\delta^c_{i-1},z))\neq\emptyset\} \right)\\
		& = \sharp \left( E_i\setminus \left \{ z<\delta_i^c;\ z\in \nu_i\left (\N_{[1,i)} \right )\right \} \right)\ \\
		& = \sharp \left(E_i\setminus  \nu_i\left (\N_{[1,i)}\right )\right),
	\end{array}
$$
where the last equality holds because $c(E_i)\leq c\left ( \nu_i\left ( \N_{[1,i)} \right )\right )$.
\fim\smallskip

In the statement of Corollary \ref{thetaf} the equality among the cardinality of the difference between the two value sets and the length of the quotient of the two fractional ideals follows from Remark \ref{r=1}.

To get neat formulas in what follows, we put $\N_{I\setminus \{1,\ldots,r\}}:=\mathcal{I}_1$. Notice that this agrees with $$\Theta_1(E)=\sharp\left ( E_1\setminus  \nu_1\left ( \N_{I\setminus \{1,\ldots,r\}}\right )\right )=\sharp (E_1\setminus\nu_1(\mathcal{I}_1))=\sharp(E_1\setminus E_1)=0.$$

As a consequence of Proposition \ref{codim}  and Corollary \ref{thetaf}, we get the following enhanced result:

\begin{cor}\label{theorem-colenth}
Let  $\mathcal{J},\mathcal{I}\subseteq\mathcal{K}$ be fractional ideals of $\mathcal{O}$ with $\mathcal{J}\subseteq\mathcal{I}$,  $D=\nu(\mathcal{J})$ and $E=\nu(\mathcal{I})$, then 
{\small \begin{equation}\label{dimensionf}
	\ell_{\mathcal{O}}\left (\frac{\mathcal{I}}{\mathcal{J}}\right )  =\sum_{i\in I}\left (\ell_i\left ( \frac{\mathcal{I}_i}{\mathcal{J}_i}\right )+\sharp \left(D_i\setminus  \nu_i\left ( \N_{[1,i)}(\mathcal{J})\right )\right)-\sharp\left( E_i\setminus  \nu_i\left ( \N_{[1,i)}(\mathcal{I})\right )\right )\right),
\end{equation}
}
or equivalently,
$$\ell_{\mathcal{O}}\left (\frac{\mathcal{I}}{\mathcal{J}}\right )  =\sum_{i\in I}\ell_i\left ( \frac{ \pi_i(\N_{[1,i)}(\mathcal{I}))}{ \pi_i(\N_{[1,i)}(\mathcal{J}))}\right ).$$
\end{cor}
\pf The expression (\ref{dimensionf}) follows directly from Proposition \ref{codim} and Corollary \ref{thetaf}.

For the second part, since $D_i=\nu_i(\mathcal{J}_i)$ and $E_i=\nu_i(\mathcal{I}_i)$, from Remark \ref{r=1}, we get 
$$\sharp \left(D_i\setminus  \nu_i\left ( \N_{[1,i)}(\mathcal{J})\right )\right)=\ell_i\left ( \frac{\mathcal{J}_i}{\pi_i\left ( \N_{[1,i)}(\mathcal{J})\right )}\right ), $$
and
$$\sharp \left(E_i\setminus  \nu_i\left ( \N_{[1,i)}(\mathcal{I})\right )\right)=\ell_i\left ( \frac{\mathcal{I}_i}{\pi_i\left ( \N_{[1,i)}(\mathcal{I})\right )}\right ).$$
In this way, 
\begin{equation}\label{conta1}
\ell_i\left ( \frac{\mathcal{I}_i}{\mathcal{J}_i}\right )+\sharp \left(D_i\setminus  \nu_i\left ( \N_{[1,i)}(\mathcal{J})\right )\right)=\ell_i\left ( \frac{\mathcal{I}_i}{\pi_i\left ( \N_{[1,i)}(\mathcal{J})\right )}\right ).
\end{equation}

Notice that $\N_J(\mathcal{J})\subseteq\N_J(\mathcal{I})$ for every $J\subset I$, so $\pi_i\left ( \N_{[1,i)}(\mathcal{J})\right )\subseteq \pi_i\left ( \N_{[1,i)}(\mathcal{I})\right )$, and consequently,
{\small
\begin{equation}\label{conta2}
\ell_i\left ( \frac{\mathcal{I}_i}{\pi_i\left ( \N_{[1,i)}(\mathcal{J})\right )}\right )-\ell_i\left ( \frac{\mathcal{I}_i}{\pi_i\left ( \N_{[1,i)}(\mathcal{I})\right )}\right )=
\ell_i\left ( \frac{ \pi_i\left (\N_{[1,i)}(\mathcal{I})\right )}{ \pi_i\left (\N_{[1,i)}(\mathcal{J})\right )}\right ).
\end{equation}
}
substituting (\ref{conta1}) and (\ref{conta2}) in (\ref{dimensionf}) concludes the proof.
\fim\medskip

We close this section with a lemma that will be useful in the sequel.

\begin{lema} \label{tec}
Let $\mathcal{O}$ be the local ring of an irreducible curve and let $\mathcal{L}\subseteq\mathcal{M}$ be fractional ideals of $\mathcal{O}$. If $h\in\mathcal{O}$ then $$\ell_{\mathcal{O}} \left ( \frac{\mathcal{M}}{h\mathcal{L}}\right )=\nu(h)+\ell_{\mathcal{O}} \left ( \frac{\mathcal{M}}{\mathcal{L}}\right )$$
where $h\mathcal{L}=\{hl;\ l\in\mathcal{L}\}\subseteq\mathcal{L}$.
\end{lema}
\pf Since $h\mathcal{L}\subseteq\mathcal{L}\subseteq\mathcal{M}$ and $\nu(h\mathcal{L})=\nu(h)+\nu(\mathcal{L})$, by Remark \ref{r=1}, we get $$\ell_{\mathcal{O}} \left ( \frac{\mathcal{M}}{h\mathcal{L}}\right )=\sharp \left ( \nu(\mathcal{M})\setminus\{\nu(h)+\nu(\mathcal{L})\}\right )=\sharp\{\gamma\in\nu(\mathcal{M});\ \gamma-\nu(h)\notin\nu(\mathcal{L})\}.$$

Since $\nu(\mathcal{L})$ has a conductor, for any $0\leq i<\nu(h)$ there exists $\gamma_i=\min\{\nu(\mathcal{L})\cap\{i+\nu(h)\mathbb{Z}\}\}$. As $\gamma_i\in\nu(\mathcal{L})\subseteq\nu(\mathcal{M})$ and $\gamma_i-\nu(h)\notin\nu(\mathcal{L})$, it follows that
$$\nu(\mathcal{M})\setminus\{\nu(h)+\nu(\mathcal{L})\}=\left (\nu(\mathcal{M})\setminus\nu(\mathcal{L})\right )\ \dot{\cup}\ \{\gamma_i;\ 0\leq i<\nu(h)\}.$$
Hence,
$$\ell_{\mathcal{O}} \left ( \frac{\mathcal{M}}{h\mathcal{L}}\right )=\nu(h)+\sharp(\nu(\mathcal{M})\setminus\nu(\mathcal{L}))=\nu(h)+\ell_{\mathcal{O}} \left ( \frac{\mathcal{M}}{\mathcal{L}}\right ).$$
\fim

\section{Tjurina number for plane curves}
In this section we will prove in a more elementary way refinements of results contained in \cite{Bayer}, in the particular case of plane curves, that we will need.

Let $\mathcal{O}=\mathbb C\{X,Y\}/\langle f \rangle$ be the local ring of an analytic reduced plane curve $C$, defined by $f=\prod_{i\in I}f_i\in\mathbb{C}\{X,Y\}$, where $I=\{1,\ldots,r\}$, $r>1$, with $f_i$ irreducible and $\langle f_i\rangle\neq\langle f_j\rangle$ for $i\neq j$. The irreducible components of $C$ defined by the $f_i$,  $i\in I$ , will be denoted by $C_i$.

The Jacobian ideal of $f$ is the ideal of $\Oh$ defined by
$$\mathcal{J}(f)=\{A f_{X}+B f_{Y};\ A,B\in\mathcal{O}\}\subseteq\mathcal{O}.$$ 

The \emph{Tjurina number} $\tau(C)$ of $C$ is the colength of $\mathcal{J}(f)$ in $\mathcal{O}$, that is, 
$$\tau(C):=\ell_{\mathcal{O}}\left ( \frac{\mathcal{O}}{\mathcal{J}(f)}\right )=\dim_{\mathbb{C}}\frac{\mathbb{C}\{X,Y\}}{\langle f, f_X, f_Y\rangle},$$
where $\ell_{\Oh}$ stands for the length of an $\Oh$-module.

Since $f_X=\sum_{i\in I}(\prod_{j\in I\setminus\{i\}}f_j)(f_i)_X$ and $f_Y=\sum_{i\in I}(\prod_{j\in I\setminus\{i\}}f_j)(f_i)_Y$, we have
\begin{equation}\label{deltaJ}
\mathcal{J}(f)=\left \{ \sum_{i\in I}\left (\prod_{j\in I\setminus\{i\}}f_j\right )\left ( A (f_i)_{X}+B (f_i)_{Y}\right );\ A,B\in\mathcal{O}\right \},
\end{equation}
and, for every $i\in I$, viewing $\J(f)$ in $\K$,
\begin{equation} \label{deltaJ_i} 
\mathcal{J}(f)_i :=\pi_i(\J(f))=\left (\prod_{j\in I\setminus\{i\}}f_j\right ) \mathcal{J}(f_i),
\end{equation}
where $\mathcal{J}(f_i)$ is the Jacobian ideal of $f_i$ in $\mathcal O_i$.

Caution: the above formula shows that the $i$-th projection of the Jacobian ideal in $\Oh$ is not the Jacobian ideal of the $i$-th branch $C_i$ of $C$. 

In the sequel
we denote by ${\rm I}(f_i,f_j)$ the intersection multiplicity ${\rm I}(C_i,C_j)$ that is,  ${\rm I}(f_i,f_j)=\dim_{\mathbb{C}}\frac{\mathbb{C}\{X,Y\}}{\langle f_i,f_j\rangle}=\nu_i(f_j)$.

\begin{rem} \label{sumthetai} In the present framework, if $\Gamma=\nu(\mathcal O)$, we have that
\[
\ell_{\mathcal{O}}\left ( \frac{\prod_{i\in I}\mathcal{O}_i}{\mathcal{O}}\right )=\sum_{i=2}^r \Theta_i(\Gamma)=\sum_{1\leq i<j\leq r}\textup{I}(f_i,f_j).
\]

Indeed, by Example \ref{prod}, Corollary \ref{thetaf} and Remark \ref{r=1}, we have that
\[ 
\begin{array}{rcl}
\ell_{\mathcal{O}}\left ( \frac{\prod_{i\in I}\mathcal{O}_i}{\mathcal{O}}\right ) &=&\sum_{i=2}^{r}\Theta_i(\Gamma)=\sum_{i=2}^{r}\sharp \left (\Gamma_i\setminus\nu_i(\N_{[1,i)}(\Oh))\right ) \\ \\
                              &=& \sum_{i=2}^{r}\sharp \left (\Gamma_i\setminus\nu_i(\prod_{1 \leq j<i}f_j)\right) =                                                
                             \sum_{i=2}^r\ell_i\left ( \frac{\mathcal{O}_i}{\left (\prod_{1\leq j<i}f_j\right )\mathcal{O}_i}\right ) \\ \\
                             &= & \sum_{i=2}^r\nu_i(\prod_{1\leq j<i}f_j)= \sum_{1\leq i< j\leq r}\textup{I}(f_i,f_j),
\end{array}
\]
where $\Gamma_i=\nu_i(\mathcal{O}_i)$, and, for $J\subset I$, we used the fact that  $\N_J(\Oh)=f^J\mathcal O$, with $f^J=\prod_{j\in J}f_j$ (see Remark \ref{Danna}). 
\end{rem} \smallskip

We have the following result:
\begin{prop} \label{tjurinaref} Let $\Delta=\nu(\mathcal J(f))$ and $\Delta_i=\nu_i(\mathcal J(f)_i)$, then one has
\[
\tau(C)=\sum_{i\in I}\tau(C_i)+\sum_{1\leq i< j\leq r}\textup{I}(f_i,f_j)+\sum_{i=2}^r\sharp \left (\Delta_i\setminus \nu_i\left (\N_{[1,i)}(\J(f))\right )\right ).
\]
\end{prop}
\pf From Proposition \ref{codim} and \eqref{deltaJ_i} one has
{\small
$$\tau(C)=\ell_{\mathcal{O}}\left (\frac{\mathcal{O}}{\mathcal{J}(f)} \right )=\sum_{i\in I}\left (\ell_i\left (\frac{\mathcal{O}_i}{\prod_{j\in I\setminus\{i\}}f_j\cdot\mathcal{J}(f_i)} \right )+\Theta_i(\Delta)-\Theta_i(\Gamma)\right ),$$
}

By Lemma \ref{tec}, we have
$$\begin{array}{ll}
\ell_i\left (\frac{\mathcal{O}_i}{\prod_{j\in I\setminus\{i\}}f_j\cdot\mathcal{J}(f_i)} \right ) & =\nu_i\left (\prod_{j\in I\setminus\{i\}}f_j\right )+\ell_i\left (\frac{\mathcal{O}_i}{\mathcal{J}(f_i)} \right )\\ \\
 & =\sum_{j\in I\setminus\{i\}} \textup{I}(f_i,f_j)+\tau(C_i),
\end{array}
$$
hence
$$\tau(C)=\sum_{i\in I}\left( \sum_{j\in I\setminus\{i\}} \textup{I}(f_i,f_j)+\tau(C_i) +\Theta_i(\Delta)-\Theta_i(\Gamma) \right).$$

Now, the result follows by Remark \ref{sumthetai} and Corollary \ref{thetaf}.
\fim\medskip

In the sequel, $J$ will denote a proper subset of $I$ and $K=I\setminus J$. Let us denote by $C_J$ and $C_K$ the plane curves determined by $f^J:=\prod_{j\in J}f_j$ and $f^K:=\prod_{k\in K}f_k$, respectively. In \cite[Corollary 12]{Bayer}, an additive formula is proved that relates $\tau(C)$, $\tau(C_J)$ and $\tau(C_K)$. Using Proposition \ref{tjurinaref} we will give an alternative way to express such an additive formula.

Let us consider the Jacobian ideal $\mathcal{J}(f^J )$ of $f^J$ in $\mathcal O_J$ and, for  $i\in J$, its $i$-th component 
\begin{equation}\label{fJi}
\mathcal{J}(f^J)_i=(\prod_{l\in J\setminus\{i\}}f_l)\cdot \mathcal{J}(f_i).
\end{equation}

For $i\in J$ and $L\subset J$, we put $\Delta^J_i:=\nu_i\left ( \mathcal{J}(f^J)\right )$ and $\N^J_L:=\N_L(\mathcal{J}(f^J))$. Then for every $i\in J$ we have 
\begin{equation}\label{conta}
\Delta^J_i=\nu_i\left ( \prod_{l\in J\setminus\{i\}}f_l\right )+\nu_i\left ( \mathcal{J}(f_i)\right ).
\end{equation}

\begin{lema}\label{lema-equal} Let $J$ be a proper subset of $I$ and $K=I\setminus J$. Then, for every $i\in J$, we have
	$$\Delta_i=\nu_i(f^K)+\Delta^J_i,\ \ \ \mbox{and}$$
	$$\nu_i(\N_{J\setminus \{i\}})=\nu_i(f^K)+\nu_i\left (\N^J_{J\setminus\{i\}}\right ).$$ In particular,
	\begin{equation}\label{equality}
		\sharp \left (\Delta_i\setminus \nu_i\left (\N_{J\setminus\{i\}}\right )\right )=\sharp \left (\Delta^J_i\setminus \nu_i\left (\N^J_{J\setminus\{i\}}\right )\right ).
		\end{equation}
\end{lema}
\pf  By (\ref{deltaJ_i}) and (\ref{conta}), for every $i\in J$, we have 
\begin{equation}\label{equal1}\begin{array}{ll}
\Delta_i & =\nu_i(\prod_{l\in I\setminus\{i\}}f_l)+\nu_i(\mathcal{J}(f_i))\vspace{0.2cm}\\
& =\nu_i(f^K)+\nu_i(\prod_{l\in J\setminus\{i\}}f_l)+\nu_i(\mathcal{J}(f_i))\vspace{0.2cm}\\ 
& =\nu_i(f^K)+\Delta^J_i.
\end{array}\end{equation} 

We proceed with the second part of the proof. If $h\in \J(f)$, from \eqref{deltaJ}, we may write
\[
h=\sum_{i\in I}\left (\prod_{j\in I\setminus{\{i\}}}f_j\right )\left (A (f_i)_{X}+B (f_i)_{Y}\right ), \ \text{with} \ A,B\in \Oh,
\]
then, because $\pi_l(f_l)=0$, we have 
\[
\pi_l(h)=\pi_l \left( \left( \prod_{j\in I\setminus \{l\}}f_j \right) (A(f_l)_X+B(f_l)_Y )\right).
\]

Now suppose that $h\in \N_l=\N_l(\J(f))$, that is, $\pi_l(h)=0$. Since $\pi_l\left( \prod_{j\in I\setminus \{l\}}f_j\right)\neq 0$, because $f_l$ does not divide $\prod_{j\in I\setminus \{l\}}f_j$, and since $\Oh_l$ is a domain, the condition $\pi_l(h)=0$ is equivalent to the condition $\pi_l(A(f_l)_X+B(f_l)_Y)=0$. Recalling that $\nu_l(h)=\nu_l(\pi_l(h))$, we conclude that 
\begin{equation}\label{equal2}
	\begin{array}{l}
\nu_i\left (\N_{J\setminus\{i\}}\right ) =\left \{
\begin{array}{l}
\nu_i(f^{I\setminus\{i\}}\left (A (f_i)_{X}+B (f_i)_{Y}\right )),\ \ \mbox{with}\\
 \pi_l(A(f_l)_{X}+B (f_l)_{Y})=0, \ \ \forall\ l\in J\setminus\{i\}
 \end{array}\right \} \vspace{0.2cm} \\

\hspace{1.3cm}  =\nu_i(f^K)+ \left \{
\begin{array}{l}
	\nu_i(f^{J\setminus\{i\}}\left (A (f_i)_{X}+B (f_i)_{Y}\right )), \ \ \mbox{with}\\
	\pi_l(A (f_l)_{X}+B (f_l)_{Y})=0,  \ \ \forall\ l\in J\setminus\{i\}
\end{array}\right \} \vspace{0.2cm} \\

\hspace{1.3cm} = \nu_i(f^K)+ \nu_i\left (\N^J_{J\setminus\{i\}}\right ).

\end{array}
\end{equation}

Equality (\ref{equality}) follows from (\ref{equal1}) and (\ref{equal2}). 
\fim\medskip

If $i\in L \subset J \subset I$, then by applying twice \eqref{equality}, we have that
\begin{equation} \label{L<J}
\sharp(\Delta_i\setminus \nu_i(\N_{L\setminus \{i\}})) = \sharp(\Delta^L_i\setminus \nu_i(\N^L_{L\setminus \{i\}}))=\sharp(\Delta^J_i\setminus \nu_i(\N^J_{L\setminus \{i\}})).
\end{equation}

Without loss of generality we may assume that $J=\{1,\ldots ,j\}$ and
$K=\{j+1,\ldots ,r\}$. We are now in position to prove the following result:

\begin{teo} \label{tjurinaIJK}
Let $C$ be an analytic reduced plane curve defined by $f=\prod_{i\in I}f_i\in\mathbb{C}\{X,Y\}$ with $I=\{1,\ldots ,r\}$. If $C_J$ and $C_K$ denote the plane curves defined by $\prod_{j\in J}f_j$ and $\prod_{k\in K}f_k$, respectively, where $J=\{1,\ldots ,j\}$ and $K=\{j+1,\ldots ,r\}$, then 
$$ \tau(C)=\tau(C_J)+\tau(C_K)+\textup{I}(C_J,C_K)+\sum_{k\in K}\sharp\left ( \nu_k\left (\N_{[j+1,k)}\right )\setminus \nu_k\left ( \N_{[1,k)}\right )\right ),$$
where $$\textup{I}(C_J,C_K)=\sum_{i\in J\atop k\in K}\textup{I}(f_i,f_k),\ \N:=\N(\mathcal{J}(f))\ \mbox{and}\ \nu_{j+1}\left (\N_{[j+1,j+1)}\right ):=\Delta_{j+1}.$$
\end{teo}
\pf From Proposition \ref{tjurinaref}, we get
$$\tau(C_J)=\sum_{i\in J}\tau(C_i)+\sum_{1\leq i< l\leq j}\textup{I}(f_i,f_l)+\sum_{i=2}^j\sharp \left (\Delta^J_i\setminus \nu_i\left (\N^J_{[1,i)}\right )\right ),\ \ \mbox{and}$$
$$\tau(C_K)=\sum_{k\in K}\tau(C_k)+\sum_{j+1\leq i< l\leq r}\textup{I}(f_i,f_l)+\sum_{i=j+2}^r\sharp \left (\Delta^K_i\setminus \nu_i\left (\N^K_{[j+1,i)}\right )\right ).$$

So, again by Proposition \ref{tjurinaref}, we have
\begin{equation}\label{aux}
\begin{array}{rl}
\tau(C)=&\tau(C_J)+\tau(C_K)+\textup{I}(C_J,C_K)+\sum_{i=2}^r\sharp \left (\Delta_i\setminus \nu_i\left (\N_{[1,i)}\right )\right ) \\ \\
        & -\sum_{i=2}^j\sharp \left (\Delta^J_i\setminus \nu_i\left (\N^J_{[1,i)}\right )\right ) -\sum_{i=j+2}^r\sharp \left (\Delta^K_i\setminus \nu_i\left (\N^K_{[j+1,i)}\right )\right ).
\end{array}
\end{equation}

By \eqref{L<J}, for $2\leq i\leq j$, by taking $L=\{1,\ldots, i\}$, we have 
\begin{equation} \label{aux1} 
\sharp \left (\Delta_i\setminus \nu_i\left (\N_{[1,i)}\right )\right )=\sharp \left (\Delta^J_i\setminus \nu_i\left (\N^J_{[1,i)}\right )\right ),
\end{equation} and for $j+2\leq i\leq r$, by taking $L=\{j+1,\ldots,i\}$, from \eqref{L<J}, we have
\begin{equation}\label{aux2}
\begin{array}{l}
\hspace{-1.75cm}\sharp \left ( \Delta_i\setminus \nu_i\left (\N_{[1,i)}\right )\right )
-\sharp \left ( \Delta^K_i\setminus \nu_i\left ( \N^K_{[j+1,i)}\right )\right )= \\ \\
\sharp \left ( \Delta_i\setminus \nu_i\left (\N_{[1,i)}\right )\right )
-\sharp \left ( \Delta_i\setminus \nu_i\left ( \N_{[j+1,i)}\right )\right )= \\ \\
\hspace{2.5cm} \sharp\left ( \nu_i\left ( \N_{[j+1,i)}\right )\setminus \nu_i\left ( \N_{[1,i)}\right )\right ).
\end{array}
\end{equation}
The result follows from \eqref{aux}, \eqref{aux1} and \eqref{aux2}.
\fim\medskip

The above theorem is a refinement of \cite[Corollary 12]{Bayer} in the case of plane curves.

\section{Proof of Dimca's conjecture}

In \cite{Dimca}, A. Dimca proposed the following conjecture:

\begin{conj}[Conjecture 1.7, \cite{Dimca}]\label{conj-dimca}
With the previous notation one has $$\tau(C)\leq\tau(C_J)+\tau(C_K)+2\cdot \textup{I}(C_J,C_K)-1.$$
\end{conj}

This section is dedicated to the proof of this conjecture. Here we will always use $\N$ to denote $\N(\J(f))$ and $J=\{1,\ldots ,j\}$. 

By Remark \ref{r=1}, for every $j+2\leq k\leq r$ we have that
$$\sharp\left ( \nu_k\left (\N_{[j+1,k)}\right )\setminus \nu_k\left (\N_{[1,k)}\right )\right ) =\ell_k\left ( \frac{\pi_k(\N_{[j+1,k)})}{\pi_k(\N_{[1,k)})}\right ).$$

Since 
$$f^J \N_{[j+1,k)}\subseteq \N_{[1,k)}\subseteq \N_{[j+1,k)},$$
by Lemma \ref{tec} one has, for every $k$, with $j+2\leq k\leq r$,
{\small \begin{equation}\label{aux3}
\sharp\left ( \nu_k\left (\N_{[j+1,k)}\right )\setminus \nu_k\left (\N_{[1,k)}\right )\right )\leq \ell_k\left ( \frac{\pi_k(\N_{[j+1,k)})}{\pi_k(f^J \N_{[j+1,k)})}\right )=\nu_k( f^J ) = \textup{I}\left ( f_k,\prod_{j\in J}f_j\right ).
\end{equation}
}

Moreover, we have the following result:

\begin{lema}\label{estimate}
With the above notations, we have
$$\sharp\left ( \Delta_{j+1}\setminus \nu_{j+1}\left (\N_{[1,j+1)}\right )\right )\leq \nu_{j+1}(f^J)-1=\textup{I}\left (f_{j+1},\prod_{j\in J}f_j\right )-1.$$
\end{lema}
\pf Initially, notice that 
\begin{equation} \label{subseteq}
f^J \mathcal{J}(f)+\left (f^J_Yf_X-f^J_Xf_Y\right )\mathcal{O}\subseteq \N_{[1,j+1)}=\N_J.
\end{equation}
In fact, for every $l\in J$, we get
\[\begin{array}{l}
\pi_l\left ( f^J \mathcal{J}(f)+\left (f^J_Yf_X-f^J_Xf_Y\right )\mathcal{O}\right )= \\ \\
\hspace{0.5cm} 0 \cdot \mathcal{J}(f)_l+ \left (\prod_{i\in J\setminus\{l\}}f_i\right )\left (\prod_{i\in I\setminus\{l\}}f_i\right )\left ( (f_l)_Y(f_l)_X-(f_l)_X(f_l)_Y\right )=0.
\end{array}
\]

From the inclusion \eqref{subseteq} and from \eqref{deltaJ_i} we get
$$\begin{array}{l}
\nu_{j+1}\left (\N_{[1,j+1)}\right )\supseteq	
\nu_{j+1}\left ( f^J \mathcal{J}(f)+\left (f^J_Yf_X-f^J_Xf_Y\right )\mathcal{O}\right )   \supseteq \\ \\
\nu_{j+1}\left ( f^J \mathcal{J}(f)\right )\bigcup \nu_{j+1}\left ( \left (f^J_Yf_X-f^J_Xf_Y\right )\mathcal{O}\right ) = \\ \\
\nu_{j+1}\left ( f^J \mathcal{J}(f)_{j+1}\right )\bigcup\ \nu_{j+1}\left ( f^{I\setminus\{j+1\}}\left (f^J_Y(f_{j+1})_X-f^J_X(f_{j+1})_Y\right )\mathcal{O}_{j+1}\right ). 
\end{array}$$

By \cite[Proposition 2.1.1]{del}, we have
\begin{equation}\label{del}
\nu_{j+1}(h_Y(f_{j+1})_X-h_X(f_{j+1})_Y)=\nu_{j+1}(h)+\mu_{j+1}-1,\ \ \ \ \forall\ h\in\langle X,Y\rangle
\end{equation}
where $\mu_{j+1}$ is the conductor of $\Gamma_{j+1}=\nu_{j+1}(\mathcal{O}_{j+1})$. Hence,
$$\begin{array}{l}
\hspace{-1.5cm}\nu_{j+1}\left ( f^{I\setminus\{j+1\}}\left (f^J_Y(f_{j+1})_X-f^J_X(f_{j+1})_Y\right )\mathcal{O}_{j+1}\right )=\vspace{0.2cm}\\
\hspace{1.5cm}=\nu_{j+1}(f^{I\setminus\{j+1\}})+\nu_{j+1}(f^J)+\mu_{j+1}-1+\Gamma_{j+1}.
\end{array}$$

In addition, by \eqref{del}, we have 
$$\nu_{j+1}\left( \langle X,Y\rangle\right) +\mu_{j+1}-1 =\Gamma_{j+1}\setminus\{0\}+\mu_{j+1}-1\subseteq\nu_{j+1}(\mathcal{J}(f_{j+1})),$$ 
and, since from \eqref{fJi},
$$\Delta_{j+1}=\nu_{j+1}\left ( \mathcal{J}(f)_{j+1}\right )=\nu_{j+1}\left ( f^{I\setminus\{j+1\}}\mathcal{J}(f_{j+1})\right ),$$  one has
$$\Gamma_{j+1}\setminus\{0\}+\nu_{j+1}\left ( f^{I\setminus\{j+1\}}\right )+\nu_{j+1}(f^J)+\mu_{j+1}-1\subseteq \nu_{j+1}\left ( f^J \mathcal{J}(f)_{j+1}\right )$$
and 
$$ \nu_{j+1}\left ( f^{I\setminus\{j+1\}}\right )+\nu_{j+1}(f^J)+\mu_{j+1}-1\not\in\nu_{j+1}\left ( f^J \mathcal{J}(f)_{j+1}\right ),$$
because being $\mu_{j+1}$ the conductor of $\Gamma_{j+1}$, we have $\mu_{j+1}-1\notin \Gamma_{j+1}$.

Consequently,
$$\begin{array}{l}
\nu_{j+1}\left ( f^J \mathcal{J}(f)_{j+1}\right )
\bigcup\ \nu_{j+1}\left ( f^{I\setminus\{j+1\}}\left (f^J_Y(f_{j+1})_X-f^J_X(f_{j+1})_Y\right )\mathcal{O}_{j+1}\right )=\vspace{0.2cm}\\
\hspace{2cm}=\nu_{j+1}\left ( f^J \mathcal{J}(f)_{j+1}\right )
\dot{\bigcup} \{\nu_{j+1}\left ( f^{I\setminus\{j+1\}}\right )+\nu_{j+1}(f^J)+\mu_{j+1}-1\}\end{array}$$
and
{\small $$\nu_{j+1}\left ( \N_{[1,j+1)}\right )\supseteq \nu_{j+1}\left ( f^J \mathcal{J}(f)_{j+1}\right )
\dot{\bigcup} \{\nu_{j+1}\left ( f^{I\setminus\{j+1\}}\right )+\nu_{j+1}(f^J)+\mu_{j+1}-1\}.$$}

Hence, by Lemma \ref{tec}, we obtain
$$\begin{array}{l}
\hspace{-1cm}\sharp\left ( \Delta_{j+1}\setminus \nu_{j+1}\left (\N_{[1,j+1)}\right )\right )\leq \vspace{0.2cm}\\
\hspace{2.5cm} \leq\sharp\left ( \Delta_{j+1}\setminus\nu_{j+1}\left ( f^J \mathcal{J}(f)_{j+1}\right )\right )-1 \vspace{0.2cm} \\
\hspace{2.5cm} =\ell_{j+1}\left ( \frac{\pi_{j+1}(\mathcal{J}(f)_{j+1})}{\pi_{j+1}(f^J \mathcal{J}(f)_{j+1})}\right )-1=\nu_{j+1}(f^J)-1.\end{array}$$ 
\fim

As a consequence of the above results, we have the following theorem that proves the Conjecture \ref{conj-dimca}.

\begin{teo}\label{proof-dimca}
	Let $C$ be an analytic reduced plane curve defined by $f=\prod_{i\in I}f_i\in\mathbb{C}\{X,Y\}$ with $I=\{1,\ldots ,r\}$. If $C_J$ and $C_K$ denote the plane curve defined by $\prod_{j\in J}f_j$ and $\prod_{k\in K}f_K$ respectively, where $J$ is a proper subset of $I$ and $K=I\setminus J$, then
	\begin{equation}\label{inequality1} \tau(C)-(\tau(C_J)+\tau(C_K))\leq 2\cdot \textup{I}(C_J,C_K)-1.\end{equation}
	In particular, we get
	$$\tau(C)-\sum_{i\in I}\tau(C_{i})\leq 2\cdot\sum_{1\leq i<j\leq r}\textup{I}(C_i,C_j)-r+1=\mu(C)-\sum_{i\in I}\mu(C_{i}),$$
	where $\mu(C)=\textup{I}(f_X,f_Y)$ and $\mu(C_i)=\textup{I}((f_i)_X,(f_i)_Y)$ denote the Milnor numbers of $C$ and $C_i$, respectively.
\end{teo}
\pf 
Considering $J=\{1,\ldots ,j\}$ and recalling that $\N_{[j+1,j+1)}=\Delta_{j+1}$, from Theorem \ref{tjurinaIJK} we have
\[\begin{array}{ll}
	\tau(C) - (\tau(C_J)+\tau(C_K)) & = \textup{I}(C_J,C_K)+\sharp\left ( \Delta_{j+1}\setminus\nu_i(\N_{[1,j+1)})\right )+ \vspace{0.2cm} \\
	&\hspace{1cm} \sum_{k=j+2}^r\sharp \left (\nu_{k}\left ( \N_{[j+1,k)}\right )\setminus \nu_k\left (\N_{[1,k)} \right )\right )\vspace{0.2cm} \\
	& \leq 2\textup{I}(C_J,C_K)-1,
	\end{array}
\]
where the last inequality follows from estimate \eqref{aux3} and Lemma \ref{estimate}.

In particular, from Proposition \ref{tjurinaref} and Lemma \ref{estimate} we get
{\small 
\[\begin{array}{ll}
\tau(C) - \sum_{i\in I}\tau(C_i) & = \sum_{1\leq i<j\leq r}\textup{I}(C_i,C_j)+\sum_{i=2}^r\sharp \left (\Delta_i\setminus \nu_i\left (\N_{[1,i)} \right )\right )\vspace{0.2cm} \\
& \leq \sum_{1\leq i<j\leq r}\textup{I}(C_i,C_j)+\sum_{i=2}^r\left ( \nu_i\left ( \prod_{1\leq j<i}f_j\right )-1\right )\vspace{0.2cm} \\
& = 2 \sum_{1\leq i<j\leq r}\textup{I}(C_i,C_j)-r+1 = \mu(C)-\sum_{i\in I}\mu(C_{i}),
\end{array}
\]}
where the last equality follows from the well known Milnor formula $\mu(C)=\sum_{i\in I}\mu(C_{i})+2 \sum_{1\leq i<j\leq r}\textup{I}(C_i,C_j)-r+1$ for plane curves (see \cite[Theorem 6.5.1]{wall}, for instance).
\fim\medskip

The bound in Theorem \ref{proof-dimca} is sharp, since, by a result of K. Saito (see \cite{saitinho}), the quasihomogeneous curve given by $(Y^n-a_1X^m)(Y^n-a_2X^m)\cdots(Y^n-a_rX^m)$, where the $a_i$ are non-zero distinct complex numbers and ${\rm GCD}(n,m)=1$, achieves the bound. So, at this point it is natural to pose the following\medskip

\noindent \textbf{Problem:} Characterize all plane curves that achieve the bound.

\section{The bridge between differentials and Jacobian ideals}

In \cite[Corollary 11]{Bayer} 
is presented a formula relating Tjurina's number of a curve $C:f=0$ to the Tjurina numbers of its components $C_i:f_i=0$ through a correcting term involving intersection multiplicities between these $C_i$'s and numerical invariants $\Theta_i$ obtained from the value $\Lambda$ set of the module of K\"ahler differentials $\Omega(\Oh)$ modulo its torsion submodule $\mathcal T$. On the other hand, we got the same formula with the same invariants by using the value set of the Jacobian ideal. This is not a surprise, since one has that 
$\Delta=c(\Gamma)-(1,\ldots,1)+\Lambda$ (cf. \cite[Corollary 6]{Bayer}). Hence, in the presence of Corollary \ref{thetaf}, it follows that Proposition \ref{tjurinaref} is a particular case of \cite[Corollary 11]{Bayer} for plane curves, and at the same time, a refinement of it.
In our new approach the $\Theta_i$ where obtained from certain submodules $\N_{J}(\J(f))$ of the Jacobian ideal $\J(f)$. The challenge now is to reinterpret these submodules in terms of the module of differentials. Establishing such connection, apparently increases the possibility of making effective computations.

For the convenience of the reader we will recall some facts from \cite[Section 2]{Bayer} and follow the same notations.

We denote the module of differentials $1$-forms of $\mathbb{C}^2$ by $\Omega^1=\mathbb{C}\{X,Y\}\textup{d}X+\mathbb{C}\{X,Y\}\textup{d}Y$ and define the submodule  $\mathcal{F}_f:=\mathbb{C}\{X,Y\}\cdot \textup{d}f+f\cdot \Omega^1$.

From \cite[(2.1)]{saito} we have the following $\Oh$-modules isomorphism:
\[
\Omega(\Oh) \simeq \dfrac{\Omega^1}{{\mathcal F}(f)}.
\]

We will consider $\Omega(\Oh)/\T$ as a fractional ideal of $\Oh$ in $\K$ as in \cite[Equation (6)]{Bayer} and one has
\[
\left( \dfrac{\Omega(\Oh)}{\T}\right)_i:=\pi_i\left( \dfrac{\Omega(\Oh)}{\T}\right)=\dfrac{\Omega(\Oh_i)}{\T_i},
\]
where $\T_i$ is the torsion submodule of $\Omega(\Oh_i)$. We will denote by $\Lambda$ the value set of the fractional ideal $\dfrac{\Omega(\Oh)}{\T}$. So, we have
\[
\Lambda_i:=\pr_i(\Lambda)=\nu_i\left( \dfrac{\Omega(\Oh_i)}{\T_i}\right).
\]

We define
\[
\Ss_i:= \dfrac{\Omega(\Oh)}{\T} \bigcap \ker \pi_i.
\]

Notice that
$$\begin{array}{rcl}
\Ss_i &=& \left\{\overline{\omega}\in \dfrac{\Omega(\Oh)}{\T}; \; \exists \overline{g} \in \Oh_i\setminus \{0\}, \ \overline{g\omega}=0\in \Omega(\Oh_i)\right\} \\  \\
&=& \left\{\overline{\omega}\in \dfrac{\Omega(\Oh)}{\T}; \; \exists g\notin \langle f_i \rangle, \ g\omega \in \F(f_i)\right\} \\
&=& \left\{ \overline{\omega}\in \dfrac{\Omega(\Oh)}{\T}; \; \omega\in\Omega^1, \ \frac{\omega\wedge \textup{d}f_i}{\textup{d}X\wedge \textup{d}Y}\in\langle f_i\rangle\right \},
\end{array}
$$
where the last equality follows from \cite[(1.1), page 266]{saito}. 

To preserve the analogy with the notations of Section 2, for $J\subset I$, we define
\[
\Ss_J=\bigcap_{j\in J} \Ss_j,
\]
and get the following result:

\begin{teo} One has 
\begin{equation} \label{1st}
\nu_i\left ( \mathcal{N}_{[1,i)}\right )=\left \{\nu_i \overline{\left (\frac{\omega\wedge\textup{d}f}{\textup{d}X\wedge\textup{d}Y}\right )};\ \overline{\omega}\in\Ss_{[1,i)} \right \}.
\end{equation}
In particular,	
$\sharp \left(\Delta_i\setminus\nu_i\left ( \mathcal{N}_{[1,i)}\right )\right)=\sharp \left(\Lambda_i\setminus\nu_i\left ( \Ss_{[1,i)}\right )\right)$.
\end{teo}
\pf We have that
\[
\begin{array}{rcl}
\mathcal{N}_j &= &\{Af_X+Bf_Y\in\mathcal{J}(f);\ \pi_j(Af_X+Bf_Y)=0\} \\ \\
&  = &\left \{\overline{\frac{\omega\wedge\textup{d}f}{\textup{d}X\wedge\textup{d}Y}}\in f_j\mathcal{O};\ \omega=B\textup{d}X-A\textup{d}Y\in\Omega^1 \right \}  \\ \\
 & =& \left \{\overline{\frac{\omega\wedge\textup{d}f}{\textup{d}X\wedge\textup{d}Y}}\in\mathcal{O};\ \overline{\omega}\in\mathcal{S}_j\right \}.
\end{array}
\]

From this, \eqref{1st} follows.\medskip

Now, since $\textup{d}f=\sum_{i\in I}\left (\prod_{j\in I\setminus{\{i\}}}f_j\right )\textup{d}f_i$, we get 
\[
\pi_i \overline{\left (\frac{\omega\wedge\textup{d}f}{\textup{d}X\wedge\textup{d}Y}\right )}=\pi_i \overline{\left (\left (\prod_{j\in I\setminus{\{i\}}}f_j\right )\cdot\frac{\omega\wedge\textup{d}f_i}{\textup{d}X\wedge\textup{d}Y}\right )}
\]
hence
\begin{equation} \label{interm}
\begin{array}{ll}
\nu_i\overline{\left (\frac{\omega\wedge\textup{d}f}{\textup{d}X\wedge\textup{d}Y}\right )} & =\nu_i\left (\prod_{j\in I\setminus{\{i\}}}f_j\right )+\nu_i\overline{\left (\frac{\omega\wedge\textup{d}f_i}{\textup{d}X\wedge\textup{d}Y}\right )}\\ \\
& =\nu_i\left (\prod_{j\in I\setminus{\{i\}}}f_j\right )+\nu_i(\omega)+\mu_i-1,
\end{array}
\end{equation}
where the last equality follows from a well known equality (see for example, \cite[Theorem 5 (iii)]{Bayer}).

From \eqref{1st} and \eqref{interm} we get 
\[
 \nu_i\left ( \N_{[1,i)}\right )=\nu_i\left ( \Ss_{[1,i)} \right )+\nu_i\left (\prod_{j\in I\setminus{\{i\}}}f_j\right )+\mu_i-1.
 \]
 
From the relation 
\[
\Delta=\Lambda+c(\Gamma)-(1,\ldots,1),
\]
we get, by projecting with $\pr_i$,
\[
\Delta_i=\Lambda_i+c(\Gamma)_i-1.
\]

Hence from \eqref{conds} we conclude the proof of the theorem.
\fim
\medskip

\textbf{Acknowledgment.}  We would like to thank Professor Alexandru Dimca who called our attention to Conjecture \ref{conj-dimca} that was the  motivation for doing this work.

\end{document}